# INVERTIBLE POLYNOMIAL MAPPINGS VIA NEWTON NON-DEGENERACY

Y. CHEN, L.R.G. DIAS, K. TAKEUCHI, AND M. TIBĂR

ABSTRACT. We prove a sufficient condition for the Jacobian problem in the setting of real, complex and mixed polynomial mappings. This follows from the study of the bifurcation locus of a mapping subject to a new Newton non-degeneracy condition.

## 1. INTRODUCTION

Unlike the local setting, the critical locus is not the only obstruction to produce diffeomorphisms in the global setting. A well-known example by Pinchuk [Pi] yields a polynomial mapping $\mathbb{R}^2 \to \mathbb{R}^2$ with no singularities but which is not invertible, thus providing a counter-example to the strong Jacobian Conjecture over the reals. A natural question, also posed by Bivià-Ausina for real mappings in [Bi], would then be: given a polynomial mapping $G : \mathbb{A}^n \to \mathbb{A}^n$ with $\text{Sing} G = \emptyset$, where $\mathbb{A} = \mathbb{R}$ or $\mathbb{C}$, under what general enough conditions $G$ becomes a diffeomorphism. For polynomial mappings over $\mathbb{C}$, being a diffeomorphism (actually "injective" is enough by [BR]) ensures that $G^{-1}$ is a polynomial map too (see [Gr, Proposition 17.9.6] and [CR]), but this fact is no more true over the reals.

We obtain here a new sufficient condition for the invertibility of $G$ as a by-product of the study of the *bifurcation locus* of a polynomial mapping $F : \mathbb{A}^n \to \mathbb{A}^k$, $n \geq k$. This is the minimal set of points $B(F) \subset \mathbb{A}^k$ such that the restriction $F_| : F^{-1}(\mathbb{A}^k \setminus B(F)) \to \mathbb{A}^k \setminus B(F)$ is a locally trivial fibration. One has no complete knowledge over this set unless $k = 1$ and $n = 2$, see [Su], [HL] [Du], [Ti1], [TZ]. In more variables, one may estimate $B(F)$ by some "reasonably good" superset $B \supset B(F)$ by using criteria of regularity at infinity. This was first done in case of complex polynomial functions $f : \mathbb{C}^n \to \mathbb{C}$ with conditions like *tameness* [Br1], *Malgrange regularity* [Pa], *ρ-regularity* [NZ2], [Ti1], *$\mathcal{W}$-equisingularity* [Ve], [ST] etc. Each of these conditions holds over $\mathbb{R}$ too and exhibits, in both settings, a finite subset of "non-regular values" containing the bifurcation set $B(F)$ (see e.g. [Ti1], [Ti2]). For $k > 1$, Rabier [Ra] considered an asymptotic regularity derived from a Palais-Smale type condition which extends the Malgrange regularity. This allows him to define a set of "asymptotically critical values" $K_\infty(F)$ and to prove that $B(F) \subset F(\text{Sing} F) \cup K_\infty(F)$. In order to evaluate $K_\infty(F)$, we introduce here a new *Newton non-degeneracy condition at infinity*, see Definition 3.4. Besides real and complex polynomial mappings, we consider *mixed mappings*, which are mappings $\mathbb{C}^n \to \mathbb{C}^k$ in

2010 *Mathematics Subject Classification.* 14D06, 58K05, 57R45, 14P10, 32S20, 58K15.
L.R.G. Dias acknowledges the Brazilian grants CAPES-Proc. 2929/10-04 and FAPESP-Proc. 2008/10563-4. M. Tibăr and L.R.G. Dias were partially supported by the USP-COFECUB grant Uc Ma 133/12.

variables $z$ and $\bar{z}$. We first obtain the following bound of the asymptotical critical locus of $F = (f_1, \ldots, f_k)$, in terms of two subsets of $\mathbb{A}^k$ depending on the Newton polyhedra of the functions $f_i$, denoted by $N(F)$ and $A(F)$ (see Definitions 3.5 and 3.8), the first of which is algebraic and the second, semi-algebraic:

**Theorem 1.1.** *Let $F : \mathbb{A}^n \to \mathbb{A}^k$, $n \geq k$, be a real, complex or mixed polynomial mapping depending effectively on all the variables, such that $F(0) = 0$. If $F$ is non-degenerate at infinity then:*
$$K_\infty(F) \subset N(F) \cup A(F).$$
□

Using the inclusions $B(F) \subset F(\operatorname{Sing} F) \cup S(F)$ and $S(F) \subset K_\infty(F)$ proved in [DRT] and discussed here in §2 and Proposition 2.4, we immediately get the following bound for the bifurcation locus:

**Corollary 1.2.** $B(F) \subset F(\operatorname{Sing} F) \cup N(F) \cup A(F)$. □

By Proposition 2.4, our Theorem 1.1 extends to mappings (real, complex or mixed) the result for complex polynomial functions proved by Némethi-Zaharia [NZ1] in the complex setting and more recently by Chen-Tibăr [CT, Theorem 1.1(a)] in the mixed setting. Moreover, when applied to mixed mappings $F : \mathbb{C}^n \to \mathbb{C}^k$, our Theorem 1.1 yields a better result than for the same underlying real mapping $F : \mathbb{R}^{2n} \to \mathbb{R}^{2k}$. This is not only trivially visible in the first term of the union $N(F) \subset \mathbb{A}^k \setminus (\mathbb{A}^*)^k$ since $\mathbb{A} = \mathbb{R}$ is replaced by $\mathbb{A} = \mathbb{C}$, but also in the second term since the involved Newton polyhedra turn out to be different. Moreover, the statement in the mixed setting cannot be deduced from the one in the real setting; its proofs is similar but needs specific notations and preliminaries. We give in §4 the proof of Theorem 1.1 and show the following significant consequence (cf. Definitions 3.1 and 3.4):

**Corollary 1.3.** *Suppose that $F$ is non-degenerate at infinity and that $f_i$ is convenient, for any $i = 1, \ldots, k$. Then $K_\infty(F) = \emptyset$.*

This provides an extension to *mappings* (which can be real, complex or mixed) of Broughton's classical result [Br1, Proposition 3.4] which tells that if a complex polynomial function $f : \mathbb{C}^n \to \mathbb{C}$ is convenient and Newton non-degenerate then $F$ is "tame", thus $K_\infty(F) = \emptyset$. The proof of our Corollary 1.3 will be given in §4.

With all these preparations, we may state and prove the announced result:

**Theorem 1.4.** *Let $F : \mathbb{A}^n \to \mathbb{A}^n$ be a $C^1$ real semi-algebraic mapping such that $\operatorname{Sing} F = \emptyset$. If $K_\infty(F) = \emptyset$ then $F$ is a global diffeomorphism.*

*In particular, if $F = (f_1, \ldots, f_n)$ is a real, or a mixed, or a complex polynomial mapping with $\operatorname{Sing} F = \emptyset$, non-degenerate at infinity, and if $f_i$ is convenient for all $i = 1, \ldots, n$, then $F$ is a global diffeomorphism.*

*Proof.* Let $J_F$ denote the set of points at which $F$ is not proper (see Definition 2.2 and [Je1, Definition 3.3]). By [KOS, Proposition 3.1], one has $K_\infty(F) = J_F$. Thus if $K_\infty(F) = \emptyset$ then $F$ is proper. It is moreover a submersion since $\operatorname{Sing} F = \emptyset$ by hypothesis. A proper submersion is an open and closed mapping, a general topological fact. It follows that $\operatorname{Im} F = \mathbb{A}^n$, therefore $F$ is a covering and it must be one-to-one since its image $\mathbb{A}^n$ is simply connected. Our first assertion follows. Remark that the final part of this proof is actually Hadamard's theorem, see e.g. [vdE, p. 240]. Our second assertion then follows by Corollary 1.3. □

One of the new issues of our paper is the non-degeneracy condition at infinity which appears to be a generic condition (Definition 3.4). This extends to mappings the definitions of "Newton non-degeneracy at infinity" for functions, both in the complex setting [Ku], [Br1], [Br2], [NZ1] and in the more recently developed mixed setting [CT]. Moreover, this works over the reals too. Our Definition 3.4 is particularly designed to treat the case of *non-convenient* polynomial mappings (cf Definition 3.1) since bifurcation values at infinity appear only in this context, see Corollary 1.3. In the real setting, Bivià-Ausina considers in [Bi] a different Newton non-degeneracy condition for real polynomial mappings $F : \mathbb{R}^n \to \mathbb{R}^k$ and uses it in case each component $f_i$ is convenient. He proves a result like the second statement of our Theorem 1.4 in case of a real polynomial mapping $F : \mathbb{R}^n \to \mathbb{R}^n$. His proof aims to obtain the properness of $F$, necessary for the bijectivity of $F$, via an interpretation of his non-degeneracy condition in terms of Łojasiewicz exponents. The author observes in [Bi, p. 746] that the techniques used in his paper work only for real polynomial mappings. We show in §5 that in the real setting and for $n = k$, the two definitions are equivalent. However, Proposition 5.2 and Example 5.3 tell that, whenever $n > k$ and $f_i$ is convenient for every $i$, our definition is strictly more general than the one considered in [Bi]. Our definition of non-degeneracy for mappings is new especially for mixed and complex mappings and in general for non-convenient mappings in all settings. The mixed setting is particularly interesting upon the real one by the fact that a mixed mapping may be convenient without the underlying real components to be convenient, whereas if all the real components are convenient then the mixed functions must be also convenient (see Example 5.4).

## 2. Asymptotic critical values

Let $F = (f_1, \ldots, f_k) : \mathbb{A}^n \to \mathbb{A}^k$, $n \geq k$, be a $C^1$ real semi-algebraic mapping with $f_i \not\equiv$ const., $\forall i = 1, \ldots, k$. The *Milnor set* of $F$, denoted by $M(F)$, is the critical locus of the mapping $(F, \rho)$, where $\rho : \mathbb{A}^n \to \mathbb{R}_{\geq 0}$ denotes the Euclidean distance function. From its definition, it follows that $M(F)$ is a closed semi-algebraic subset of $\mathbb{A}^n$ and that, in the case $n = k$, $M(F)$ coincides with the whole $\mathbb{A}^n$.

**Definition 2.1.** Let $F : \mathbb{A}^n \to \mathbb{A}^k$ be a $C^1$ non-constant real semi-algebraic mapping. The set of *asymptotic non $\rho$-regular values* of $F$ is defined as

$$S(F) := \left\{ c \in \mathbb{A}^k \mid \exists \{x_l\}_{l \in \mathbb{N}} \subset M(F), \lim_{l \to \infty} \|x_l\| = \infty \text{ and } \lim_{l \to \infty} F(x_l) = c \right\}.$$

**Definition 2.2.** [Je1, Definition 3.3] Let $F : \mathbb{A}^n \to \mathbb{A}^k$ be a continuous mapping. We say that $F$ is proper at a point $c \in \mathbb{A}^k$ if there exists an open neighbourhood $U$ of c such that the restriction $F_{|F^{-1}(U)} : F^{-1}(U) \to U$ is a proper mapping. We denote by $J_F$ the set of points at which $F$ is not proper.

One also has the *set of generalised critical values* $K(F) := F(\mathrm{Sing}\, F) \cup K_\infty(F)$, where:

**Definition 2.3.** [Ra, p. 670], [KOS, p. 68] The set of *asymptotic critical values* of a real semi-algebraic mapping $F : \mathbb{A}^n \to \mathbb{A}^k$ is defined as:

$$K_\infty(F) := \{c \in \mathbb{A}^k \mid \exists \{x_l\}_{l \in \mathbb{N}} \subset \mathbb{A}^n, \lim_{l \to \infty} \|x_l\| = \infty,$$
$$\lim_{l \to \infty} F(x_l) = c \text{ and } \lim_{l \to \infty} (1 + \|x_l\|) \nu(dF(x_l)) = 0\},$$

where $\nu(B) := \inf_{\|\varphi\|=1} \|B^*(\varphi)\|$, for $B \in \mathcal{L}(\mathbb{A}^n, \mathbb{A}^k)$.

In the context of semi-algebraic $F \colon \mathbb{R}^n \to \mathbb{R}^p$, Kurdyka, Orro and Simon showed in [KOS] that $K_\infty(F)$ is a semi-algebraic set of dimension at most $k - 1$. Gaffney [Ga] defined a generalised Malgrange condition in the setting of complex polynomial mappings $\mathbb{C}^n \to \mathbb{C}^k$ and proved that this condition yields a set $A_G(F)$ of non-regular values such that $B(F) \subset F(\operatorname{Sing} F) \cup A_G(F)$. Then Jelonek [Je2] showed that the asymptotic conditions employed in [KOS] and in [Ga] are equivalent, i.e. that $K_\infty(F) = A_G(F)$.

More recently, [DRT] proved, in the setting of semi-algebraic mappings, the inclusion $B(F) \subset F(\operatorname{Sing} F) \cup S(F)$ and that $S(F)$ and $F(\operatorname{Sing} F) \cup S(F)$ are closed semi-algebraic sets of dimension $\leq k - 1$ (more precisely, [DRT, Theorem 5.7] in the case $n > k$ and [KOS, Theorem 3.1 and Proposition 3.1] in the case $n = k$).

From the above definitions we get the inclusions $S(F) \subset J_F$ and $K_\infty(F) \subset J_F$. In the case $n = k$ one has $J_F = K_\infty(F)$ by [KOS, Proposition 3.1], hence $S(F) = J_F = K_\infty(F)$.

For $n \geq k$, the inclusion $S(F) \subset K_\infty(F)$ was shown in [CT, Proposition 2.2] for a mixed polynomial, and in the more general real setting in [DRT, Corollary 5.8]. Here we offer a new and direct proof of this inclusion.

**Proposition 2.4.** *Let $F = (f_1, \ldots, f_k) : \mathbb{A}^n \to \mathbb{A}^k$ be a $C^1$ non-constant real semi-algebraic mapping with $n \geq k$. Then $S(F) \subseteq K_\infty(F)$.*

*Proof.* We give the proof over $\mathbb{R}$. Then the statement over $\mathbb{C}$ can be obtained from the one over $\mathbb{R}$ by using the identification $\mathbb{C}^n$ with $\mathbb{R}^{2n}$.

In the case $n = k$, as explained just above, we have equality. We concentrate in the following to the case $n > k$.

Let $c = (c_1, \ldots, c_k) \in S(F)$. Since $M(F)$ is semi-algebraic, one can use the Curve Selection Lemma at infinity to find an analytic path $\phi = (\phi_1, \ldots, \phi_n) : ]0, \epsilon[ \to M(F) \subset \mathbb{R}^n$ such that $\lim_{t \to 0} \|\phi(t)\| = \infty$ and $\lim_{t \to 0} F(\phi(t)) = c$.

Since $\phi(t) \in M(F)$ if and only if $\operatorname{rank} dF(\phi(t)) < k+1$, there exist curves $\lambda(t), b_1(t), \ldots, b_k(t)$ such that $(\lambda(t), b_1(t), \ldots, b_k(t)) \neq (0, \ldots, 0), \forall t$, and one has the equality:

$$(1) \qquad \lambda(t)(\phi_1(t), \ldots, \phi_n(t)) = b_1(t) \frac{\partial f_1}{\partial x}(\phi(t)) + \ldots + b_k(t) \frac{\partial f_k}{\partial x}(\phi(t)),$$

where $\frac{\partial f_i}{\partial x}(\phi(t)) = \left( \frac{\partial f_i}{\partial x_1}(\phi(t)), \ldots, \frac{\partial f_i}{\partial x_n}(\phi(t)) \right)$, for $i = 1, \ldots, k$.

Consider $b(t) = (b_1(t), \ldots, b_k(t))$. From the equality (1) and by the statements that $(\lambda(t), b_1(t), \ldots, b_k(t)) \neq (0, \ldots, 0), \forall t$, and $\lim_{t \to 0} \|\phi(t)\| = \infty$, we have $b(t) \neq 0, \forall t$ and consequently, from (1), we obtain:

$$(2) \qquad \frac{\lambda(t)}{\|b(t)\|}(\phi_1(t), \ldots, \phi_n(t)) = \frac{b_1(t)}{\|b(t)\|} \frac{\partial f_1}{\partial x}(\phi(t)) + \ldots + \frac{b_k(t)}{\|b(t)\|} \frac{\partial f_k}{\partial x}(\phi(t));$$

and we will denote $\lambda_0(t) := \frac{\lambda(t)}{\|b(t)\|}$ and $a(t) := \frac{b(t)}{\|b(t)\|}$. So, $\|a(t)\| = 1$ and one obtains the following equalities:

$$(3) \qquad \sum_{i=1}^{k} a_i(t)\frac{\mathrm{d}}{\mathrm{d}t} f_i(\phi(t)) = \left\langle \sum_{i=1}^{k} a_i(t)\frac{\partial f_i}{\partial x}(\phi(t)), \phi'(t) \right\rangle = \frac{1}{2}\lambda_0(t)\frac{\mathrm{d}}{\mathrm{d}t}\|\phi(t)\|^2,$$

where the later follows from (2), i.e., from the equality $\sum_{i=1}^{k} a_i(t)\frac{\partial f_i}{\partial x}(\phi(t)) = \lambda_0(t)\phi(t)$.

On the other hand, since $\lim_{t\to 0} f_i(\phi(t)) = c_i$, it follows that $\mathrm{ord}_t\left(\frac{\mathrm{d}}{\mathrm{d}t}f_i(\phi(t))\right) \geq 0$, $i = 1, \ldots, k$. This and the equality (3) imply:

$$(4) \qquad 0 \leq \mathrm{ord}_t\left(\lambda_0(t)\frac{\mathrm{d}}{\mathrm{d}t}\|\phi(t)\|^2\right) < \mathrm{ord}_t\left(\lambda_0(t)\|\phi(t)\|^2\right).$$

Now, from (2) one obtains:

$$(5) \qquad \mathrm{ord}_t\left(\|\phi(t)\|\|a_1(t)\frac{\partial f_1}{\partial x}(\phi(t)) + \ldots + a_k(t)\frac{\partial f_k}{\partial x}(\phi(t))\|\right) = \mathrm{ord}_t\left(|\lambda_0(t)|\|\phi(t)\|^2\right),$$

which is positive by (4). This implies:

$$\lim_{t\to 0} \|\phi(t)\|\|a_1(t)\frac{\partial f_1}{\partial x}(\phi(t)) + \ldots + a_k(t)\frac{\partial f_k}{\partial x}(\phi(t))\| = 0,$$

which, in turn, implies $\lim_{t\to 0} \|\phi(t)\|\nu(\mathrm{d}F(\phi(t))) = 0$. This shows that $\mathrm{c} \in K_\infty(F)$. $\square$

3. NEWTON POLYHEDRA AND THE NON-DEGENERACY AT INFINITY

Let $f : \mathbb{A}^n \to \mathbb{A}$ be a non-constant polynomial function, $\mathbb{A} = \mathbb{R}$ or $\mathbb{C}$. We write:

$$f(\mathrm{x}) = \sum_\nu c_\nu \mathrm{x}^\nu,$$

where $\nu = (\nu_1, \cdots, \nu_n) \in \mathbb{N}^n$ and $\mathrm{x}^\nu = x_1^{\nu_1} \cdots x_n^{\nu_n}$. In the mixed case, we write:

$$f(\mathrm{z}, \bar{\mathrm{z}}) = \sum_{\nu,\mu} c_{\nu,\mu} \mathrm{z}^\nu \bar{\mathrm{z}}^\mu,$$

where $\nu = (\nu_1, \cdots, \nu_n), \mu = (\mu_1, \cdots, \mu_n) \in \mathbb{N}^n, \mathrm{z}^\nu = z_1^{\nu_1} \cdots z_n^{\nu_n}$ and $\bar{\mathrm{z}}^\mu = \bar{z}_1^{\mu_1} \cdots \bar{z}_n^{\mu_n}$.

Mixed polynomials is a much larger class than complex polynomials and have been introduced by Oka, who studied several aspects of their local topology in [Oka2], [Oka3] and some other more recent articles. In [CT] we have used the *mixed Newton polyhedron at infinity* of a polynomial function.

**Definition 3.1.** Let $f : \mathbb{A}^n \to \mathbb{A}$ be a non-constant polynomial function (resp. mixed function) such that $f(0) = 0$. We call $\mathrm{supp}(f) = \{\nu \in \mathbb{N}^n \mid c_\nu \neq 0\}$ (resp. $\mathrm{supp}(f) = \{\nu + \mu \in \mathbb{N}^n \mid c_{\nu,\mu} \neq 0\}$) the *support* of $f$. We say that $f$ is *convenient* if the intersection of $\mathrm{supp}(f)$ with each coordinate axis is non-empty. We denote by $\overline{\mathrm{supp}(f)}$ the convex hull of the set $\mathrm{supp}(f)$. The *Newton polyhedron* of $f$, denoted by $\Gamma_0(f)$, is the convex hull of the set $\{0\} \cup \mathrm{supp}(f)$. The *Newton boundary at infinity* of $f$, denoted by $\Gamma^+(f)$, is the union of the faces of the polyhedron $\Gamma_0(f)$ which do not contain the origin. By "face" we mean face of any dimension. Let $\Delta$ be such a face of $\overline{\mathrm{supp}(f)}$. The restriction of $f$ to $\Delta \cap \mathrm{supp}(f)$, denoted by $f_\Delta$, is defined as follows $f_\Delta(\mathrm{x}) := \sum_{\nu \in \Delta \cap \mathrm{supp}(f)} c_\nu \mathrm{x}^\nu$ (resp. $f_\Delta(\mathrm{z}, \bar{\mathrm{z}}) := \sum_{\nu + \mu \in \Delta \cap \mathrm{supp}(f)} c_{\nu,\mu} \mathrm{z}^\nu \bar{\mathrm{z}}^\mu$).

Let us consider in the following a real, mixed or complex mapping $F = (f_1, \ldots, f_k) : \mathbb{A}^n \to \mathbb{A}^k$, $n \geq k$ with $F(0) = 0$.

**Definition 3.2.** For some vector $\mathbf{p} = (p_1, \ldots, p_n) \in \mathbb{Z}^n$ with $p := \min_{1 \leq i \leq n} p_i < 0$, let $l_\mathbf{p}(\mathrm{v}) := \sum_{i=1}^n p_i v_i$ be the linear form defined by $\mathbf{p}$. Let then $\Delta_\mathbf{p}^j$ be the maximal face of $\Gamma_0(f_j)$ (maximal with respect to the inclusion of faces) where $l_\mathbf{p}(\mathrm{v})$ takes its minimal value on $\Gamma_0(f_j)$. We consider the following equivalence relation on the set of vectors $\mathbf{p}$ as above:
$$\mathbf{p}' \sim \mathbf{p} \iff \Delta_{\mathbf{p}'}^j = \Delta_\mathbf{p}^j, \text{ for all } 1 \leq j \leq k.$$

The equivalent classes yield a partition of $\mathbb{R}^n \setminus \mathbb{R}_+^n$ into finitely many locally closed (and not necessarily convex) cones. We may however subdivide each equivalence class into a finite number of convex polyhedral cones.

Let $\mathcal{C}(F)$ denote the finite set of cones obtained in this way. We call it *the dual subdivision associated to $F$*.

It follows from the above definition that, for any $j$, the face $\Delta_\mathbf{p}^j$ of $\Gamma_0(f_j)$ is independent of the defining vector $\mathbf{p}$ in its equivalence class, thus we may use the notations $\Delta_\sigma^j$ for $\sigma \in \mathcal{C}(F)$, instead of $\Delta_\mathbf{p}^j$ for some $\mathbf{p} \in \sigma$. The following sets are therefore well defined:
$$I_\sigma = \{1 \leq j \leq k \mid 0 \notin \Delta_\sigma^j\}, \qquad J_\sigma = \{1 \leq j \leq k \mid \Delta_\sigma^j = \{0\}\}.$$

REMARK 3.3. Let $d_\mathbf{p}^j \in \mathbb{Z}$ denote the minimum value of the restriction of $l_\mathbf{p}$ to $\overline{\mathrm{supp}(f)}$. We have the following relations which follow directly from the definitions:
  (a) $j \in I_\sigma \iff d_\mathbf{p}^j < 0$ for any $\mathbf{p} \in \sigma \iff \Delta_\sigma^j$ is a face of $\Gamma^+(f_j)$.
  (b) $j \in J_\sigma \iff \Gamma_0(f_j) \setminus \{0\}$ is included into the positive half-space defined by the hyperplane $\{l_\mathbf{p} = 0\}$ in $\mathbb{R}^n$ with normal vector $\mathbf{p}$, for any $\mathbf{p} \in \sigma$.

The following definition of non-degeneracy is inspired from Oka's work [Oka1] on complex local complete intersections and from the definition used by Matsui-Takeuchi [MT] and Esterov-Takeuchi [ET] in the global setting of complex polynomials. It was proved in [Oka1] that, in the complex context, this is a generic condition.

**Definition 3.4.** We say that the polynomial mapping $F = (f_1, \ldots, f_k) \colon \mathbb{A}^n \longrightarrow \mathbb{A}^k$ is *non-degenerate at infinity* if, for any $\sigma \in \mathcal{C}(F)$ such that $I_\sigma \neq \emptyset$ and $J_\sigma = \emptyset$, the subvariety:
$$(6) \qquad G_\sigma = \{x \in (\mathbb{A}^*)^n \mid f_{\Delta_\sigma^j}(x) = 0 \text{ for any } j \in I_\sigma\}$$
of $(\mathbb{A}^*)^n$ is a *non-degenerate complete intersection* (i.e. $\mathrm{Sing}\,(f_{\Delta_\sigma^j})_{j \in I_\sigma} \cap G_\sigma = \emptyset$), where $f_{\Delta_\sigma^j}$ is a short notation for the restriction of $f_j$ to $\Delta_\sigma^j$.

**Definition 3.5.** Let $\mathcal{C}(F)_{ex} := \{\sigma \in \mathcal{C}(F) \mid J_\sigma \neq \emptyset\}$ be called *the set of exceptional cones*. Let $N_\sigma = \{(z_1, \ldots, z_k) \in \mathbb{A}^k \mid z_j = 0 \text{ for any } j \in J_\sigma\}$. We then define the following algebraic subset of $\mathbb{A}^k$ of $\mathbb{A}$-codimension $\geq 1$:
$$N(F) := \cup_{\sigma \in \mathcal{C}(F)_{ex}} N_\sigma.$$

REMARK 3.6. The above definition implies that we have the inclusion $N(F) \subset \mathbb{A}^k \setminus (\mathbb{A}^*)^k$. Whenever $k \geq 2$, one can characterise the situations when this inclusion is strict, as follows. Let $C_i$ denote the set of hyperplanes $H \subset \mathbb{R}^n$ through the origin such that $H \cap \Gamma_0(f_i) = \{0\}$. Then: $N(F) \subsetneq \mathbb{A}^k \setminus (\mathbb{A}^*)^k$ if and only if there exists $j \in \{1, \ldots, k\}$ such that $C_j \subset \cup_{i \neq j} C_i$. Indeed, the condition $C_j \subset \cup_{i \neq j} C_i$ is equivalent to the fact that there is no $\sigma \in \mathcal{C}(F)$ such that $J_\sigma = \{j\}$.

REMARK 3.7. From the definition of $N(F)$ we immediately get the equivalence:
$N(F) = \emptyset \iff$ for any $i \in \{1, \ldots, k\}$, $f_i$ is convenient.
If all the cones $\mathbb{R}_+ \Gamma_0(f_i)$ coincide but are different from $(\mathbb{R}_{\geq 0})^n$, then $N(F) = \{0\}$. In particular, in case $k = 1$ we get $N(f) = \{0\}$ for any non-convenient polynomial $f$. Let us remark that the set $\{0\}$ appears as a component in the union of sets which occur as bound for the bifurcation set of a polynomial map $B(f)$ in the formula by Némethi-Zaharia [NZ1] and also in the one by Chen-Tibăr [CT].

Consider now cones $\sigma \in \mathcal{C}(F) \setminus \mathcal{C}(F)_{ex}$ and let $I_\sigma^c := \{1, 2, \ldots, k\} \setminus I_\sigma$.

**Definition 3.8.** Let $\mathcal{C}(F)_{aty} := \{\sigma \in \mathcal{C}(F) \setminus \mathcal{C}(F)_{ex} \mid I_\sigma^c \neq \emptyset\}$ be called *the set of atypical cones*. For some ordered set $J := \{j_1, \ldots, j_r\} \subset \{1, \ldots, k\}$, let $\pi_J : \mathbb{A}^k \to \mathbb{A}^{|J|}$ denote the projection $(x_1, \ldots, x_k) \mapsto (x_{j_1}, \ldots, x_{j_r})$. We consider the following restriction of $F$:

$$F_\sigma := (f_{\Delta_\sigma^j})_{j \in I_\sigma^c} : G_\sigma \longrightarrow \mathbb{A}^{|I_\sigma^c|},$$

its discriminant set $\operatorname{Disc} F_\sigma \subset \mathbb{A}^{|I_\sigma^c|}$ and its inverse image $A_\sigma := \pi_{I_\sigma^c}^{-1}(\operatorname{Disc} F_\sigma) \subset \mathbb{A}^k$. We then define the following semi-algebraic subset of $\mathbb{A}^k$, of $\mathbb{A}$-codimension $\geq 1$:

$$A(F) := \cup_{\sigma \in \mathcal{C}(F)_{aty}} A_\sigma.$$

REMARK 3.9. In case $k = 1$, $F = f$, one has the notion of "bad faces" of $\overline{\operatorname{supp} f}$ in [NZ1] and [CT]. Let us then remark that the "bad faces" are among the faces $\Delta_\sigma \cap \overline{\operatorname{supp} f}$ for $\sigma \in \mathcal{C}(f)_{aty}$, and that the $\sigma \in \mathcal{C}(F)_{aty}$ such that $\Delta_\sigma \cap \overline{\operatorname{supp} f}$ is not a "bad face" yields $A_\sigma = \emptyset$.

**Definition 3.10.** We say that *$F$ depends effectively on all the variables*, if for every variable $z_i$ there exists some $j(i) \in \{1, \ldots, k\}$ such that $f_{j(i)}$ depends effectively on $z_i$. This condition is natural since if it is not satisfied then our polynomial map depends on less than $n$ variables.

## 4. Proof of Theorem 1.1 and some consequences

For $I \subset \{1, \ldots, n\}$, we define $\mathbb{A}^I := \{z = (z_1, \ldots, z_n) \mid z_i = 0, i \notin I\}$, $(\mathbb{A}^*)^I := \{z = (z_1, \ldots, z_n) \mid z_i = 0 \iff i \notin I\}$ and $F^I := F_{|\mathbb{A}^I}$, the restriction of $F$ on $\mathbb{A}^I$.

The proof will be given in the mixed setting only, since the proof in the real setting follows faithfully the same pattern and only needs adapted notations.

Let $c = (c_1, \ldots, c_k) \in K_\infty(F) \setminus N(F)$. We may apply the *Curve Selection Lemma* at infinity (see [NZ1] and [CT]), namely there exists an analytic path $z(t) = (z_1(t), \ldots, z_n(t))$ defined on a small enough interval $]0, \varepsilon[$, such that $\lim_{t \to 0} \|z(t)\| = \infty$, $\lim_{t \to 0} F(z(t), \bar{z}(t)) = c$ and

(7) $$\lim_{t \to 0} \|z(t)\| \|\nu(dF(z(t)))\| = 0.$$

We have $F: \mathbb{R}^{2n} \to \mathbb{R}^{2k}$ where $z_j = x_j + iy_j$, $f_l = g_l + ih_l$. By the proof of [CT, Lemma 2.1] one has:

$$(a_i + ib_i)\frac{\partial \overline{f_l}}{\partial \overline{z}_j} + (a_i - ib_i)\frac{\partial f_l}{\partial \overline{z}_j} = a_i\frac{\partial g_l}{\partial x_j} + b_i\frac{\partial h_l}{\partial x_j} + i(a_i\frac{\partial g_l}{\partial y_j} + b_i\frac{\partial h_l}{\partial y_j}),$$

which shows that: $\nu(\mathrm{d}F(\mathrm{z})) = \inf \|\sum_{i=1}^{k}(\mu_i \overline{\mathrm{d}f_i}(\mathrm{z},\overline{\mathrm{z}}) + \overline{\mu_i}\overline{\mathrm{d}}f_i(\mathrm{z},\overline{\mathrm{z}}))\|$, for $\mu_i \in \mathbb{C}$ with $\sum_{i=1}^{k}|\mu_i|^2 = 1$. Therefore (7) yields, for any $i \in \{1,\ldots,n\}$:

(8) $\lim_{t \to 0} \|z_i(t)\| \|\mu_1(t)\frac{\overline{\partial f_1}}{\partial z_i}(\mathrm{z}(t),\overline{\mathrm{z}}(t)) + \overline{\mu_1}(t)\frac{\partial f_1}{\partial \overline{z}_i}(\mathrm{z}(t),\overline{\mathrm{z}}(t)) + \cdots + \overline{\mu_k}(t)\frac{\partial f_k}{\partial \overline{z}_i}(\mathrm{z}(t),\overline{\mathrm{z}}(t))\| = 0,$

where $\mu_j(t) \in \mathbb{C}$ and $\sum_{j=1}^{k}|\mu_j(t)|^2 = 1$, since the left hand side of (8) is less than or equal to $\|\mathrm{z}(t)\|\|\nu(\mathrm{d}F(\mathrm{z}(t)))\|$. Let $L = \{l \in \{1,\ldots,n\} \mid z_l(t) \not\equiv 0\}$. Observe that $L \neq \emptyset$ since $\lim_{t \to 0}\|\mathrm{z}(t)\| = \infty$, and write:

(9) $\quad\quad z_l(t) = z_l t^{p_l} + \text{h.o.t.}, \quad \text{where } z_l \in \mathbb{C}^*, \ p_l \in \mathbb{Z}, \ \forall l \in L.$

Consider the expansion of $F(\mathrm{z}(t),\overline{\mathrm{z}}(t))$ for all $i = 1,\ldots,k$, we have either:

$$f_i(\mathrm{z}(t),\overline{\mathrm{z}}(t)) \equiv c_i$$

or

(10) $\quad\quad f_i(\mathrm{z}(t),\overline{\mathrm{z}}(t)) = c_i + \text{h.o.t.}.$

One may assume (eventually after a change of coordinates) that $L = \{1,\ldots,m\}$ and $p = p_1 \leq p_2 \leq \cdots \leq p_m$. Notice that, since $\lim_{t \to 0}\|\mathrm{z}(t)\| = \infty$, one has $p = \min_{i \in L}\{p_i\} < 0$, which was an assumed condition in the definition of the set $\mathcal{C}(F)$ in the preceding section. Let $\mathrm{z}_0 := (z_1,\ldots,z_m,0,\ldots,0) \in (\mathbb{C}^*)^L$ and consider the linear function $l_{\mathbf{p}}(v) = \sum_{i=1}^{m}p_i v_i + \sum_{j=m+1}^{n}gv_j$, where $\mathbf{p} := (p_1,\ldots,p_m,g,\ldots,g) \in \mathbb{Z}^n$ with $g \in \mathbb{N}$ big enough. Let $\Delta_{\mathbf{p}}^{iL}$ be the *maximal face* of $\overline{\mathrm{supp}(f_i^L)}$ where $l_{\mathbf{p}}$ restricted to $\overline{\mathrm{supp}(f_i^L)}$ takes its minimal value, which we denote by $d_{\mathbf{p}}^{iL}$. We observe that, by definition of the vector $\mathbf{p}$ and by definition of $f_i^L$, one has $\Delta_{\mathbf{p}}^{iL} = \Delta_{\mathbf{p}}^{i}$, $d_{\mathbf{p}}^{iL} = d_{\mathbf{p}}^{i}$, and consequently $f_{i\Delta_{\mathbf{p}}^{iL}}^{L} = f_{i\Delta_{\mathbf{p}}^{i}}$ (in fact, for any $(v_1,\ldots,v_n) \in \overline{\mathrm{supp}(f_i)} \setminus \overline{\mathrm{supp}(f_i^L)}$, the value of $\sum_{i=1}^{m}p_i v_i + g\sum_{i=m+1}^{n}v_i$ is greater than $d_{\mathbf{p}}^{iL}$, $\forall i = 1,\ldots,k$). So we may denote $\Delta_{\mathbf{p}}^{iL}$ (resp. $d_{\mathbf{p}}^{iL}$) only by $\Delta_{\mathbf{p}}^{i}$ (resp. $d_{\mathbf{p}}^{i}$).

We have:

(11) $\quad\quad f_i(\mathrm{z}(t),\overline{\mathrm{z}}(t)) = f_i^L(\mathrm{z}(t),\overline{\mathrm{z}}(t)) = f_{\Delta_{\mathbf{p}}^{i}}^{L}(\mathrm{z}_0,\overline{\mathrm{z}}_0) t^{d_{\mathbf{p}}^{i}} + \text{h.o.t.}$

Since $\lim_{t \to 0}F(\mathrm{z}(t),\overline{\mathrm{z}}(t)) = \mathrm{c} \in \mathbb{C}^k \setminus N(F)$, one has $d_{\mathbf{p}}^{i} \leq 0$ for all $i = 1,\ldots,k$. We write:

(12) $\quad\quad \mu_i(t) = \mu_i t^{q_i} + \text{h.o.t.}, \quad \text{where } \mu_i \in \mathbb{C}^* \text{ and } q_i \geq 0.$

If $\mu_i \equiv 0$, we put $q_i = \infty$ in (12). Let $I = \left\{i \in \{1,\ldots,k\} \mid q_i + d_{\mathbf{p}}^{i} = \min_{1 \leq i \leq k}(q_i + d_{\mathbf{p}}^{i})\right\}$. As $\sum_{i=1}^{k}\|\mu_i(t)\|^2 = 1$, we have $\min_{1 \leq i \leq k}q_i = 0$. Hence $I \neq \emptyset$ and $\mu_i(t) \not\equiv 0$ for $i \in I$. We

therefore conclude $q_i + d_{\mathbf{p}}^i \leq 0$, $\forall i \in I$. Then (8) yields, for any $l \in L$:

$$\sum_{i \in I}(\mu_i z_l \overline{\frac{\partial f_{\Delta_{\mathbf{p}}^i}^L}{\partial z_l}}(\mathbf{z}_0, \overline{\mathbf{z}}_0) + \overline{\mu}_i z_l \frac{\partial f_{\Delta_{\mathbf{p}}^i}^L}{\partial \overline{z}_l})(\mathbf{z}_0, \overline{\mathbf{z}}_0) t^{q_i + d_{\mathbf{p}}^i} + \text{h.o.t.} \to 0. \tag{13}$$

Comparing the orders of the two sides in the above formula, we obtain, for any $l \in L$:

$$\sum_{i \in I}(\mu_i \overline{\frac{\partial f_{\Delta_{\mathbf{p}}^i}^L}{\partial z_l}}(\mathbf{z}_0, \overline{\mathbf{z}}_0) + \overline{\mu}_i \frac{\partial f_{\Delta_{\mathbf{p}}^i}^L}{\partial \overline{z}_l}(\mathbf{z}_0, \overline{\mathbf{z}}_0)) = 0. \tag{14}$$

Let $\mathbf{z}_1 := (z_1, \ldots, z_m, 1, \ldots, 1)$. By (14), by the definitions of the vectors $\mathbf{z}_0$ and $\mathbf{z}_1$, and by the equality $f_{\Delta_{\mathbf{p}}^i}^L = f_{\Delta_{\mathbf{p}}^i}$ explained in the paragraph before equation (11), one concludes that $\mathbf{z}_1 \in \operatorname{Sing} F_{\Delta_{\mathbf{p}}} \cap (\mathbb{C}^*)^n$. Notice that we have the equivalence: $d_{\mathbf{p}}^j < 0 \Leftrightarrow j \in I_\sigma$ for $\sigma \ni \mathbf{p}$. By (10) and (11), and since $d_{\mathbf{p}}^j < 0$, we must have $f_{\Delta_{\mathbf{p}}^j}(\mathbf{z}_1, \overline{\mathbf{z}}_1) = 0$. Therefore $\mathbf{z}_1$ belongs to the set $\operatorname{Sing}(f_{\Delta_{\mathbf{p}}^j})_{j=1}^k \cap \{\mathbf{z} \in (\mathbb{C}^*)^n \mid f_{\Delta_{\mathbf{p}}^j}(\mathbf{z}, \overline{\mathbf{z}}) = 0, \forall j \in I_\sigma\}$ and notice that this set is equal to $\operatorname{Sing} F_\sigma$ defined in the preceding section, by the non-degeneracy Definition 3.4 and since $J_\sigma = \emptyset$. If $I_\sigma^c = \emptyset$ then $\operatorname{Sing} F_\sigma = \emptyset$ by the same non-degeneracy condition, and if $I_\sigma^c \neq \emptyset$ then $\operatorname{Disc} F_\sigma$ contributes to the set $A(F)$ of Definition 3.8. Indeed, whenever $J_\sigma = \emptyset$, we have the equivalence $i \in I_\sigma^c \Leftrightarrow d_{\mathbf{p}}^i = 0$ for $\sigma \ni \mathbf{p}$, and therefore $c_i = f_{\Delta_\sigma^i}(\mathbf{z}_1, \overline{\mathbf{z}}_1) \in f_{\Delta_\sigma^i}(\operatorname{Sing} F_\sigma)$. This completes our proof. $\square$

REMARK 4.1. In [CDT] we have proved an inclusion similar to the one in Theorem 1.1 but for a different definition of non-degeneracy at infinity. The first difference is that instead of the set $N(F)$ above there occurs the larger set $\mathbb{A}^n \setminus (\mathbb{A}^*)^n$, and the inclusion $N(F) \subset \mathbb{A}^n \setminus (\mathbb{A}^*)^n$ might be strict, see Remarks 3.6 and 3.7. The second difference is between the non-degeneracy conditions. Our above definition of non-degeneracy concerns a much smaller number of faces than that in [CDT]. It yields genericity, therefore it is more natural.

Here we get however a larger number of "atypical faces" and therefore a larger set $A(F)$ than that obtained in [CDT]. However it is of the same nature and has the same minimal codimension 1.

Another similar result, obtained recently by Nguyen in the complex setting only [Ng], appears to be weaker than the one in [CDT] since it uses an even stronger definition of non-degeneracy at infinity than [CDT].

4.1. **Proof of Corollary 1.3.** Since $f_i$ is convenient for every $i$, we have that $N(F) = \emptyset$ (see Remark 3.7), and moreover $\mathcal{C}(F)_{aty} = \emptyset$, which implies $A(F) = \emptyset$. Then our statement follows from Theorem 1.1. $\square$

## 5. NON-DEGENERACY AT INFINITY AND EXAMPLES

The non-degeneracy condition formulated in the real setting by Bivià-Ausina's [Bi, Definition 3.5] appears to be equivalent to the following:

**Definition 5.1.** The mapping $F : \mathbb{R}^n \to \mathbb{R}^k$ is *non-degenerate at infinity* if the following condition is satisfied for any $\mathbf{p} = (p_1, \ldots, p_n) \in \mathbb{Z}^n$ such that $p = \min_{1 \le i \le n} p_i < 0$:

$$\tag{15} \left\{ \mathrm{x} \in (\mathbb{R}^*)^n \mid f_{\Delta_{\mathbf{p}}^j}(\mathrm{x}) = 0, \text{ for all } j = 1, \ldots, k \right\} = \emptyset.$$

Indeed, in our constructions we have used the minimal value of the linear function $l_{\mathbf{p}}(v) = \sum_{i=1}^n p_i v_i$ on $\overline{\mathrm{supp}(f_i)}$, since we have considered analytic curves depending on $t \to 0$, while in [Bi] the author used the maximal value of the linear function $l_{\mathbf{p}}(v) = \sum_{i=1}^n p_i v_i$ on $\overline{\mathrm{supp}(f_i)}$ since he considered analytic curves of variable $t \to \infty$. Modulo this difference, the original definition in [Bi] coincides to the above.

Let us first prove the relations between our non-degeneracy condition and Bivià-Ausina's. Next we give several examples illustrating the fact that our definition applies to a larger class of mappings.

**Proposition 5.2.** *Suppose that $F = (f_1, \ldots, f_k) : \mathbb{R}^n \to \mathbb{R}^k$, $k \le n$ is a polynomial mapping and that $f_i$ is convenient, for all $i = 1, \ldots, k$. If $F$ is non-degenerate at infinity after Bivià-Ausina's definition Definition 5.1, then it is also non-degenerate at infinity after our Definition 3.4.*

*This becomes an equivalence whenever $k = n$.*

*Proof.* We use the notations of the previous section. Let us fix a vector $\mathbf{p} = (p_1, \ldots, p_n) \in \mathbb{Z}^n$ with $p = \min_{1 \le i \le n} p_i < 0$. Since $f_i$ is convenient for any $i = 1, \ldots, k$, the minimal value $d_{\mathbf{p}}^i$ of $l_{\mathbf{p}}(v)$ must be strictly negative on $\overline{\mathrm{supp}(f_i)}$ and therefore $J_\sigma = \emptyset$ and $I_\sigma = \{1, \ldots, k\}$ for $\sigma \ni \mathbf{p}$, in other words $I_\sigma^c = \emptyset$. The first conclusion of our proposition follows by comparing the condition (15) with the non-degeneracy condition Definition 3.4.

In the case $n = k$, let us assume that $F$ is degenerate at infinity under Definition 5.1. Then there exists $\mathrm{x} \in (\mathbb{R}^*)^n$ and a vector $\mathbf{p} = (p_1, \ldots, p_n) \in \mathbb{Z}^n \setminus \{0\}$ with $\min_{1 \le i \le n} p_i < 0$ such that $f_{\Delta_{\mathbf{p}}^i}(\mathrm{x}) = 0$, for every $i$. This means that $I_\sigma^c = \emptyset$ and that $\mathrm{x} \in G_\sigma$ for $\sigma \ni \mathbf{p}$.

By the Euler relation for weighted-homogeneous functions, we have $\langle \mathrm{d} f_{\Delta_\sigma^i}(\mathrm{x}), \mathbf{p}\mathrm{x}\rangle = d_{\mathbf{p}}^i f_{\Delta_\sigma^i}(\mathrm{x}) = 0$, for $i = 1, \ldots, n$, where $\mathbf{p}\mathrm{x} := (p_1 x_1, \ldots, p_n x_n) \ne 0$. This implies that $(\mathrm{d} f_{\Delta_{\mathbf{p}}})_{i=1}^k(\mathrm{x}) \, \mathbf{p}\mathrm{x} = 0$, which yields $\mathrm{x} \in \mathrm{Sing}\,(\mathrm{d} f_{\Delta_{\mathbf{p}}})_{i=1}^k$. Thus $F$ is degenerate after Definition 3.4. Together with the first statement of our Proposition, this establishes the equivalence of the two definitions in the case $n = k$. □

The following example shows that the first implication in above proposition is not an equivalence in general.

EXAMPLE 5.3. Let $f \colon \mathbb{A}^2 \to \mathbb{A}$, $f(x, y) = x^2 - y^2$. Then $f$ is convenient, $f$ is non-degenerate at infinity after Definition 3.4, but degenerate after Definition 5.1.

Let us give two types of examples where our Theorem 1.4 applies beyond the results in [Bi].

EXAMPLE 5.4. Let $G = (G_1, G_2) \colon \mathbb{C}^2 \to \mathbb{C}^2$, $G(z_1, z_2, \overline{z}_1, \overline{z}_2) = (z_1 + \overline{z}_2, z_1 - \overline{z}_2)$, be a mixed polynomial mapping. Then $\mathrm{Sing}\,G = \emptyset$, $G$ is non-degenerate at infinity, and $G_1, G_2$ are convenient. Hence Theorem 1.4 applies showing that $G$ a diffeomorphism (which is already clear since $G$ is invertible). However, if one considers the associated real mapping of $G$, one has $G(x_1, x_2, y_1, y_2) = (x_1+x_2, y_1-y_2, x_1-x_2, y_1+y_2)$ and consequently $f_1, \ldots, f_4$ are not convenient and therefore Bivià-Ausina's results do not apply in this situation.

EXAMPLE 5.5. Let $F = (f_1, f_2, f_3) \colon \mathbb{A}^3 \to \mathbb{A}^3$ be a polynomial mapping defined by $F(x, y, z) = (x + yz + xy^2, y, xy + z)$, where $\mathbb{A} = \mathbb{R}$ or $\mathbb{A} = \mathbb{C}$. None of the functions $f_i$ is convenient. However we find $\operatorname{Sing} F = \emptyset$ and $K_\infty(F) = \emptyset$ by direct computations, inspite the fact that $N(F)$ is not empty. The first part of our Theorem 1.4 applies to show that $F$ is a diffeomorphism, whereas Bivià-Ausina's results do not apply in such a situation.

ICMC, Universidade de São Paulo, Av. Trabalhador São-Carlense, 400 - CP Box 668, 13560-970 São Carlos, São Paulo, Brazil and Laboratoire Paul Painlevé, Université Lille 1, 59655 Villeneuve d'Ascq, France.

*E-mail address*: yingchen@icmc.usp.br

ICMC, Universidade de São Paulo, Av. Trabalhador São-Carlense, 400 - CP Box 668, 13560-970 São Carlos, São Paulo, Brazil and Laboratoire Painlevé, UMR 8524 CNRS, Université de Lille 1, 59655 Villeneuve d'Ascq, France.

*E-mail address*: lrgdias@icmc.usp.br

Institute of Mathematics, University of Tsukuba, 1-1-1, Tennodai, Tsukuba, Ibaraki, 305-8571, Japan.

*E-mail address*: takemicro@nifty.com

Mathématiques, Laboratoire Paul Painlevé, Université Lille 1, 59655 Villeneuve d'Ascq, France.

*E-mail address*: tibar@math.univ-lille1.fr